
\input amssym.def
\input amssym
\input epsf
\magnification=1100
\baselineskip = 0.25truein
\lineskiplimit = 0.01truein
\lineskip = 0.01truein
\vsize = 8.5truein
\voffset = 0.2truein
\parskip = 0.10truein
\parindent = 0.3truein
\settabs 12 \columns
\hsize = 5.4truein
\hoffset = 0.4truein

\setbox\strutbox=\hbox{%
\vrule height .708\baselineskip
depth .292\baselineskip
width 0pt}
\font\caps=cmcsc10

\font\bigtenrm=cmr10 at 14pt

\def\sqr#1#2{{\vcenter{\vbox{\hrule height.#2pt
\hbox{\vrule width.#2pt height#1pt \kern#1pt
\vrule width.#2pt}
\hrule height.#2pt}}}}
\def\square{\mathchoice\sqr46\sqr46\sqr{3.1}6\sqr{2.3}4}

\centerline{\bigtenrm ADDING HIGH POWERED RELATIONS}
\centerline{\bigtenrm TO LARGE GROUPS}
\tenrm
\vskip 14pt
\centerline{MARC LACKENBY}
\vskip 18pt

\tenrm
\centerline{\caps 1. Introduction}
\vskip 6pt

It is a common theme in group theory that if a group $G$ is quotiented by
a sufficiently high power of an element $g \in G$, then properties of
$G$ are inherited by the quotient group. A major example of this phenomenon
is the theorem of Gromov ([4], see also [2]) asserting that if $G$ is torsion-free and
word-hyperbolic, then $G/\langle \! \langle g^n \rangle \! \rangle$
is also word-hyperbolic for all sufficiently big integers $n$. In this paper,
we will prove that this phenomenon occurs for the property
of being large. Recall that $G$ is {\sl large} if some finite index
subgroup admits a surjective homomorphism onto a non-abelian free group.
Large groups have many interesting and useful properties, including
super-exponential subgroup growth and infinite virtual first Betti
number. They are also particularly important in low-dimensional
topology. Our main theorem is the following.

\noindent {\bf Theorem 1.1.} {\sl Let $G$ be a finitely generated, large
group and let $g_1, \dots, g_r$ be a collection of elements of $G$. Then
for infinitely many integers $n$, $G/\langle \! \langle g_1^n, \dots, g_r^n \rangle \! \rangle$
is also large. Indeed, this is true when $n$ is any sufficiently large multiple of $[G:H]$,
where $H$ is any finite index normal subgroup of $G$ that admits a surjective
homomorphism onto a non-abelian free group.}

There are obvious counter-examples if `infinitely many integers $n$' is replaced
by `all but finitely many integers $n$' in the above statement. For example,
let $G = {\Bbb Z}_m \ast {\Bbb Z}_m$, where $m > 2$ and ${\Bbb Z}_m$ is the
cyclic group of order $m$. Let $g$ be a generator
of one of the free factors. Then $G$ is large, but $G/\langle \! \langle g^n \rangle \! \rangle$
is finite cyclic if $n$ and $m$ are coprime. However, when $G$ is free,
then we can obtain such a result.

\noindent {\bf Theorem 1.2.} {\sl Let $F$ be a finitely generated, free
non-abelian group. Let $g_1, \dots, g_r$ be a collection of elements of
$F$. Then, for all but finitely many integers $n$,
$F/\langle \! \langle g_1^n, \dots, g_r^n \rangle \! \rangle$ is large.
}

Theorem 1.1 is in fact a fairly rapid consequence
of Theorem 1.2. Like many results about large groups, Theorem 1.2 has
a topological proof. One realises $F$ as the fundamental group of
a bouquet of circles, and the quotient group $F/\langle \! \langle g_1^n, \dots, g_r^n 
\rangle \! \rangle$ as the fundamental group of the 2-complex $K$
obtained by attaching 2-cells along loops representing
$g_1^n, \dots, g_r^n$. The largeness of $F/\langle \! \langle g_1^n, \dots,
g_r^n \rangle \! \rangle$ is proved by analysing finite-sheeted covering
spaces of $K$ corresponding to finite index subgroups of 
$F/\langle \! \langle g_1^n, \dots, g_r^n \rangle \! \rangle$.
Specifically, the following result of the author (Theorem 3.7 of [6]) will be crucial.

\noindent {\bf Theorem 1.3.} {\sl Let $K$ be a finite connected cell complex, and let
$A$ and $B$ be subcomplexes such that $K = A \cup B$. Let $p$ be a prime and let ${\Bbb F}_p$
be the field of order $p$. Suppose that both of the maps
$$\eqalign{
H^1(A; {\Bbb F}_p) & \rightarrow H^1(A \cap B; {\Bbb F}_p) \cr
H^1(B; {\Bbb F}_p) & \rightarrow H^1(A \cap B; {\Bbb F}_p) \cr}$$
induced by inclusion are not injections. In the case $p =2$,
suppose also that the kernel of at least one of these maps
has dimension more than one. Then $\pi_1(K)$ is large.}

Theorem 1.1 is related to
the following result of the author about 3-manifolds (Theorem 3.7 of [5]).

\noindent {\bf Theorem 1.4.} {\sl Let $K$ be a non-trivial knot in
the 3-sphere, and let $m$ be any integer more than 2. Then, for all
sufficiently large integers $n$, the $mn$-fold cyclic cover of $S^3$ 
branched over $K$ has large fundamental group.}

Theorem 1.4 does not follow from Theorem 1.1, but Theorem 1.1 is
powerful enough to show that the $n$-fold cyclic cover of
$S^3$ branched over $K$ has large fundamental group, for
infinitely many positive integers $n$. The deduction runs as follows. It is
a result of Cooper, Long and Reid that $\pi_1(S^3 - K)$
is large (Theorem 1.3 of [1]). Let $M_n$ be the $n$-fold cyclic cover of $S^3$ branched
over $K$. Then $\pi_1(M_n)$ is an index $n$ subgroup of 
$\pi_1(S^3 - K)/ \langle \! \langle \mu^n \rangle \! \rangle$,
where $\mu$ is a representative for a meridian for $K$. Thus, Theorem 1.1 implies
that $\pi_1(M_n)$ is large for infinitely many $n$. 

However, Theorem 1.1 does have some interesting applications to
3-manifold theory, which go beyond Theorem 1.4. In Section 3,
we use it to prove that certain 3-manifolds constructed via Dehn
surgery have large fundamental group. In fact, we will
prove a result (Theorem 3.1) which has the following corollary.
(See Section 3 for an explanation of the relevant terminology.)

\noindent {\bf Theorem 1.5.} {\sl Let $M$ be a compact orientable
3-manifold with boundary a collection of tori. Suppose that
for one collection of slopes $(s_1, \dots, s_r)$,
with one $s_i$ on each component of $\partial M$,
$M(s_1, \dots, s_r)$ has large fundamental group.
Then this is true for infinitely many distinct collections
$(s_1, \dots, s_r)$.}

Theorem 3.1 provides much more precise information about which 
slopes $(s_1, \dots, s_r)$ yield 3-manifolds with large fundamental
group.

A result along these lines was proved by Dunfield and Thurston
(Theorem 7.3 in [3]).
They considered the case where $\partial M$ is a single torus,
and where Dehn filling along one slope $r$ yields a Seifert
fibre space $\Sigma$ with hyperbolic base orbifold. They also needed to
make a hypothesis about how $\pi_1(\partial M)$ maps into $\pi_1(\Sigma)$.
Using the largeness of $\pi_1(\Sigma)$, they proved that
for infinitely many other slopes $s$ on $\partial M$, $\pi_1(M(s))$
is large. In fact, they showed that this is true when
$\Delta(r,s)$ is big enough. Thus, when $M$ has two such
slopes $r$, $\pi_1(M(s))$ is large for all but finitely
many slopes $s$. A notable example of such an $M$ is the
exterior of the figure-eight knot. Although there are some
similarities between their approach and ours, there are also
some key differences. Our techniques only provide a proof
of a weak form of their result. However, their arguments required
detailed knowledge of the fundamental groups of Seifert fibre
spaces, and so do not readily apply when $M$ has no Seifert
fibred Dehn fillings.

\vskip 18pt
\centerline{\caps 2. Proofs of the main theorems}
\vskip 6pt

We start by explaining why Theorem 1.1 is a consequence of Theorem 1.2.
Let $G$ be a finitely generated group that is large. Let $H$ be a finite
index subgroup that admits a surjective homomorphism $\psi$ onto a non-abelian
free group $F$. Let $g_1, \dots, g_r$ be a collection of elements of $G$.
Our aim is to show that $G/\langle \! \langle g_1^n, \dots, g_r^n \rangle \! \rangle$
is large for infinitely many positive integers $n$.

Let $k_1, \dots, k_s$
be a set of representatives for the right cosets of $H$ in $G$.
Let $m$ be the smallest positive integer such that
$k_j g_i^m k_j^{-1}$ lies in $H$, for each $i$ and $j$.
Note that $m$ is finite, and when $H$ is normal in $G$, $m$ is the
lowest common multiple of the orders of $g_i H$ in $G/H$.
In particular, $m$ divides the index $[G:H]$ in
this case. Let $n$ be any positive integer, and
let $G_{mn} = \langle \! \langle g_1^{mn}, \dots, g_r^{mn} \rangle \! \rangle$
be the subgroup of $G$ normally generated by $g_1^{mn}, \dots, g_r^{mn}$.
Note that this lies in $H$, and is in fact the subgroup
of $H$ normally generated by $\{ k_j g_i^{mn} k_j^{-1} : 
1 \leq i \leq r, 1 \leq j \leq s \}$. Now,
$\{ \psi(k_j g_i^{m} k_j^{-1}) : 1 \leq i \leq r, 1 \leq j \leq s \}$
is a collection of elements of $F$. Hence,
by Theorem 1.2, 
$$F/ \langle \! \langle 
\psi(k_j g_i^{mn} k_j^{-1}) : 1 \leq i \leq r, 1 \leq j \leq s 
\rangle \! \rangle$$
is large, for all but finitely many positive integers $n$.
This is a homomorphic image of $H/G_{mn}$, which is
a finite index subgroup of $G/G_{mn}$. Thus,
$G/G_{mn} = G/ \langle \! \langle g_1^{mn}, \dots, g_r^{mn} \rangle \! \rangle$
is also large when $n$ is sufficiently big.

Note that this proof also gives some information about the set of
values of $n$ for which 
$G / \langle \! \langle g_1^{n}, \dots, g_r^{n} \rangle \! \rangle$
is large. It includes all sufficiently big multiples
of $m$, where $m$ is the positive integer defined above.
When $H$ is normal in $G$, we have already seen that
$m$ divides $[G:H]$. This proves Theorem 1.1. 

We now embark on the proof of Theorem 1.2. Let $F$ be a finitely generated,
free non-abelian group. Let $g_1, \dots, g_r$ be the given elements of $F$.
Our aim is to show that $F/\langle \! \langle g^n_1, \dots, g^n_r \rangle \! \rangle$
is large for all but finitely many positive integers $n$.

Let $\phi \colon F \rightarrow {\Bbb Z}$
be projection onto the first free factor, and let $\phi_n \colon F \rightarrow
{\Bbb Z}_n$ be the composition of $\phi$ with reduction modulo $n$. 
Realise $F$ as the fundamental group of a bouquet of circles $X$.
Let $X_n$ be the $n$-fold cyclic covering space of $X$ corresponding
to the kernel of $\phi_n$. (See Figure 1.)

\vskip 18pt
\centerline{
\epsfxsize=1.7in
\epsfbox{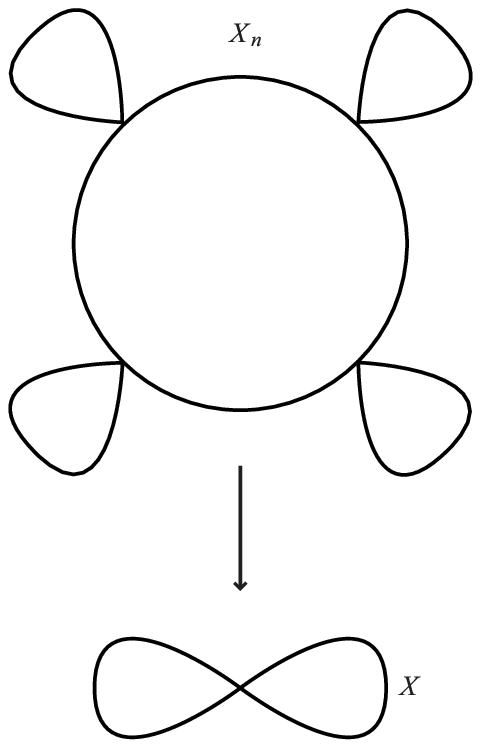}
}
\vskip 6pt
\centerline{Figure 1.}

The elements $g_1, \dots, g_r$ can be expressed as reduced words 
in $F$. They therefore determine loops in $X$, to which
we give the same names. Let $K_n$
be the 2-complex obtained by attaching 2-cells to $X$ along the words
$g_1^n, \dots, g_r^n$. Then the surjective homomorphism
$\phi_n \colon F \rightarrow {\Bbb Z}_n$ induces a surjective
homomorphism $\pi_1(K_n) \rightarrow {\Bbb Z}_n$. Let $\tilde K_n$
be the covering space of $K_n$ corresponding to the kernel of this
latter homomorphism.
This is a cell complex obtained by attaching 2-cells to $X_n$
along every lift of $g_i^n$, as $i$ runs from 1 to $r$.

Some 2-cells of $\tilde K_n$ have the same boundary loops, in the following
situation. Suppose that $\phi_n(g_i)$ is non-zero in ${\Bbb Z}_n$.
Then the word $g_i$ is an arc in $X_n$ that starts and ends at different points. 
In fact, viewing the vertices of $X_n$ as indexed by
${\Bbb Z}_n$, the start and end point of $g_i$ differ by $\phi_n(g_i)$.
So, the loop spelling $g_i^n$ starting at $0$ and the loop spelling 
$g_i^n$ starting at $\phi_n(g_i)$ have exactly the
same itineraries, but with different starting points. The same is true for 
the lift of $g_i^n$ starting at
any element of $\langle \phi_n(g_i) \rangle$, where $\langle \phi_n(g_i) \rangle$
denotes the subgroup of ${\Bbb Z}_n$ generated by $\phi_n(g_i)$. Thus, we have
$|\langle \phi_n(g_i) \rangle|$ 2-cells all with the same boundary
loop. Clearly, we could discard all but one of these 2-cells
without changing the fundamental group. Thus, we partition the
2-cells of $\tilde K_n$ into collections, where two cells belong
to the same collection if and only if their attaching loops
are lifts of the same $g_i^n$ and their starting points differ
by an element of $\langle \phi_n(g_i) \rangle$. Let $\overline K_n$
be the 2-complex obtained from $X_n$ by attaching just one
2-cell from each collection.
This has the same fundamental group as $\tilde K_n$.
Thus, in $\overline K_n$, the number of 2-cells attached along
lifts of $g_i^n$ is $[{\Bbb Z}_n : \langle \phi_n(g_i) \rangle]$.
This is equal to
$$\cases{(n, |\phi(g_i)|) & if $\phi(g_i) \not= 0$; \cr
n & if $\phi(g_i) = 0$,}$$
where $(n, |\phi(g_i)|)$ is the highest common factor
of $n$ and $|\phi(g_i)|$. Note that when $\phi(g_i) \not= 0$,
this is at most $|\phi(g_i)|$.

In the argument that follows, we will treat the elements $g_i$
such that $\phi(g_i) = 0$ in a different manner from those where
$\phi(g_i) \not= 0$. Call the former elements \break {\sl {$\phi$-trivial}},
and the latter {\sl $\phi$-non-trivial}. (See Figures 2 and 3.)
If the attaching loop of a 2-cell of $\overline K_n$ is a lift of $g_i^n$ and
$\phi(g_i) = 0$, then we also describe the 2-cell as
{\sl $\phi$-trivial}; otherwise it is {\sl $\phi$-non-trivial}.
Note that, by the above calculation, the number of
$\phi$-non-trivial 2-cells in $\overline K_n$ has an
upper bound that is independent of $n$.

\vskip 18pt
\centerline{
\epsfbox{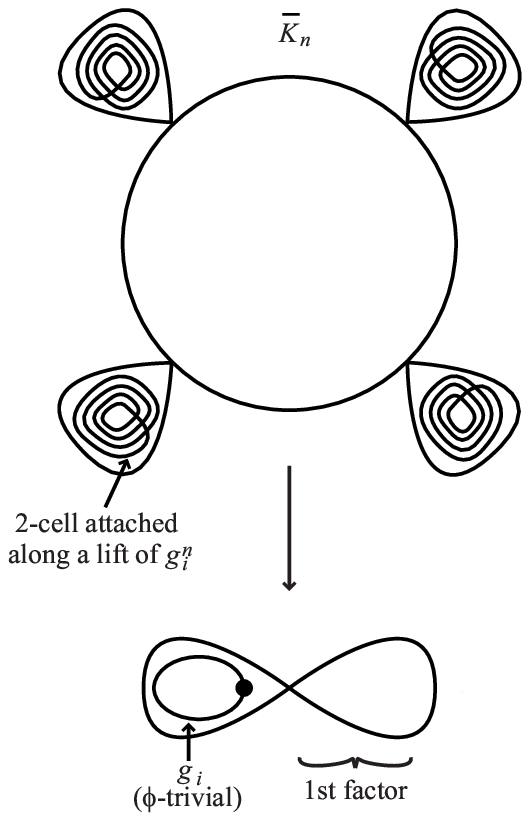}
}
\vskip 6pt
\centerline{Figure 2.}

\vskip 18pt
\centerline{
\epsfbox{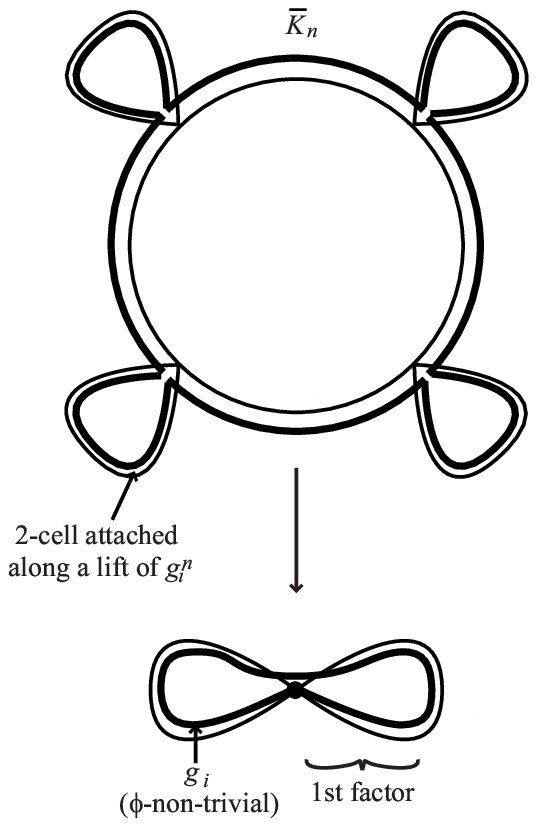}
}
\vskip 6pt
\centerline{Figure 3.}
\vfill\eject

Let $p$ be a prime number that divides $n$.
The plan is to decompose $\overline K_n$ into two subsets
$A_n$ and $B_n$ satisfying the conditions in Theorem 1.3.
(They will not in fact be subcomplexes of $\overline K_n$,
but this can easily be arranged by subdividing the cell
structure of $\overline K_n$.) Thus, an application
of Theorem 1.3 will give that $\pi_1(\overline K_n)$ is
large, and hence so is
$F/\langle \! \langle g_1^n, \dots, g_r^n \rangle \! \rangle$.

We first specify the intersection of $A_n$ and $B_n$ with
$X_n$, the 1-skeleton of $\overline K_n$. 
Note that $X_n$ is obtained from a circular graph with $n$
edges and $n$ vertices by attaching a collection of loops
to each vertex. We may label the $n$ edges of the circle
with the integers modulo $n$. Place the dividing line
between $A_n$ and $B_n$ at the midpoints of the edges
labelled $0$ and $\lceil n / 2 \rceil$. Let $P$ be
the union of these two points. Specify $A_n$ to contain
the edges increasing from $1$ to $\lceil n / 2 \rceil -1$,
together with half of the edges labelled $0$ and $\lceil n / 2 \rceil$.
In addition, if any vertex lies in $A_n$, then so do the loops
that are attached to it.
Let $B_n \cap X_n$ be the remainder of $X_n$. (See Figure 4.)

\vskip 18pt
\centerline{
\epsfbox{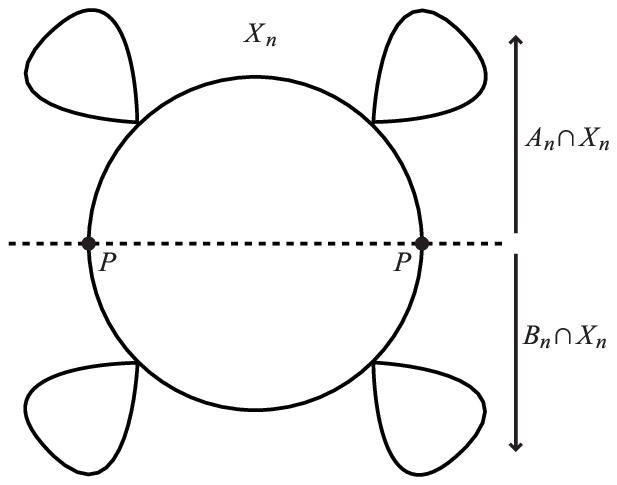}
}
\vskip 6pt
\centerline{Figure 4.}

We now define how $A_n$ and $B_n$ lie within the 2-cells of
$\overline K_n$. If the boundary of such a 2-cell $C$ lies
entirely within $A_n$ (respectively, $B_n$), then place that
entire 2-cell in $A_n$ (respectively, $B_n$). If not, place
a point of $A_n \cap B_n$ at the centre of $C$. Connect this
to each point of $A_n \cap B_n$ on the boundary of $C$ via a
radial arc. Now define $A_n$ (respectively, $B_n$) on $C$, by
defining it as the cone on $\partial C \cap A_n$ (respectively,
$\partial C \cap B_n$). (See Figure 5.)

\vskip 18pt
\centerline{
\epsfbox{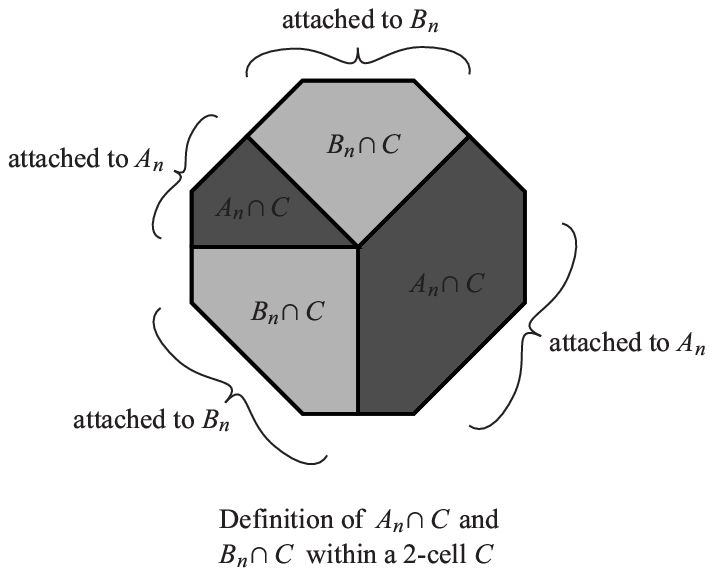}
}
\vskip 6pt
\centerline{Figure 5.}

In what follows, it will be useful to define a quantity $\Delta_i$
associated with each $g_i$. Consider the infinite cyclic
cover $X_\infty$ of $X$ corresponding to the kernel
of $\phi$. Now, $\phi$ defines
a function from the 0-skeleton of $X_\infty$ onto ${\Bbb Z}$.
Extending this linearly over each 1-cell of $X_\infty$, this
gives a function $X_\infty \rightarrow {\Bbb R}$ which we also call $\phi$.
Let $\Delta_i$ be the difference between the maximum and minimum
values of $\phi$ on a lift of $g_i$. This is clearly independent
of the choice of lift of $g_i$. (See Figure 6.)

\vskip 18pt
\centerline{
\epsfbox{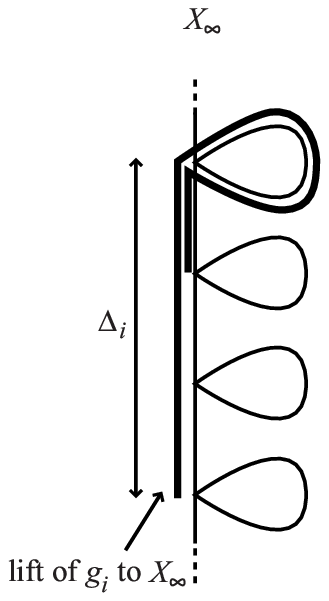}
}
\vskip 6pt
\centerline{Figure 6.}

\noindent {\sl Claim 1.} There is an upper bound,
independent of $n$, on the number of times the
$\phi$-non-trivial 2-cells of $\overline K_n$ run over
$P$.

Each $\phi$-non-trivial 2-cell of $\overline K_n$ is attached
along a lift of $g_i^n$, where $\phi(g_i) \not= 0$.
Let $j$ be the integer modulo $n$ corresponding to the
starting point of this lift. For any integer $k$
between $0$ and $n$, $g_i^k$ lifts to a sub-arc of this
loop, starting at $j$. Its endpoint is $j + k \phi_n(g_i)$, modulo $n$.
Now, for any given integer modulo $n$, $j + k \phi_n(g_i)$ 
equals that integer for at most $(n, |\phi(g_i)|)$
values of $k$. This is at most $|\phi(g_i)|$, which
is independent of $n$. Hence, there is an upper bound, independent of
$n$, on the number of integers $k$ modulo $n$ for which
$|j+ k \phi_n(g_i)| \leq \Delta_i +1$ or 
$|j+ k \phi_n(g_i) - \lceil n/2 \rceil |\leq \Delta_i +1$. But when
neither of these inequalities is satisfied, the sub-arc
of $g_i^n$ between $g_i^k$ and $g_i^{k+1}$ cannot
run over $P$. Thus, there is an upper bound,
independent of $n$, on the number of times that this lift of $g_i^n$
can run over $P$. Since there is
a universal upper bound on the number of 
$\phi$-non-trivial 2-cells in $\overline K_n$,
the claim is proved.

We now specify a graph $\Gamma_n$ lying in $\overline K_n$.
It is defined as the union of two sets:
\item{(i)} the intersection of $A_n \cap B_n$
with the $\phi$-non-trivial 2-cells, and
\item{(ii)} the 1-cells of $X_n$ lying in 
the boundary of each $\phi$-trivial 2-cell that
intersects $P$.

\noindent{\sl Claim 2.} There is an upper bound
on $d_p(\Gamma_n)$, independent of $n$.

Here, $d_p(\ \cdot \ )$ denotes the dimension of
$H_1(\ \cdot \ ; {\Bbb F}_p)$.
Note that $\Gamma_n$ has the structure of a graph.
Hence, the number of edges of $\Gamma_n$ forms
an upper bound on $d_p(\Gamma_n)$. Its
edges come in two types: radial arcs in $\phi$-non-trivial
2-cells running up to $P$, and edges in the
boundary of $\phi$-trivial 2-cells that intersect $P$.
By Claim 1, there is an upper bound, independent
of $n$, on the number of edges of the first type. 
Each $\phi$-trivial 2-cell runs $n$ times over a lift of $g_i$.
So, every edge of $\Gamma_n$ of the second type lies within
a distance $\max_i \Delta_i$ of $P$. There is therefore
a uniform upper bound (independent of $n$) on the number of
such edges. Thus, the claim is proved.

In the following claim, we show that $d_p(\overline K_n)$
is bounded below by a linear function of $n$.

\noindent {\sl Claim 3.} The following inequality holds:
$$d_p(\overline K_n) \geq n(d(F) - 1) + 1 - \sum_i{ |\phi(g_i)|}.$$

Note that $X_n$ is a connected graph with $n$ vertices and $n d(F)$ edges.
Hence, $d_p(X_n) = 1 - \chi(X_n) = n(d(F) - 1) + 1$. Now, $\overline K_n$
is obtained from $X_n$ by attaching 2-cells along lifts of $g_i^n$.
When $g_i$ is $\phi$-trivial, each corresponding 2-cell runs around a
loop $n$ times. So, attaching these does not affect $d_p$. 
(Recall that $p$ divides $n$.) When $g_i$
is $\phi$-non-trivial, we attach $(n, |\phi(g_i)|)$ 2-cells. This is
at most $|\phi(g_i)|$, and so the claim is proved.

\noindent {\sl Claim 4.} There is a constant $c$, independent of
$n$, such that the kernel of
$$H^1(\overline K_n; {\Bbb F}_p) \rightarrow H^1(A_n \cap B_n;
{\Bbb F}_p)$$
has dimension at least $n(d(F) -1) + 1- c$.

Here, and throughout the rest of the paper, homomorphisms between
cohomology groups will be induced by inclusion, when the context is clear.
By Claims 2 and 3, there is a constant $c$, independent 
of $n$, such that the kernel of
$$H^1(\overline K_n; {\Bbb F}_p) \rightarrow 
H^1(\Gamma_n; {\Bbb F}_p)$$
has dimension at least $n(d(F) - 1) + 1 - c$.
We will show that this kernel lies in the kernel of
$$H^1(\overline K_n; {\Bbb F}_p) \rightarrow H^1(A_n \cap B_n;
{\Bbb F}_p).$$ 
This will prove the claim.
Let $D$ be the union of the $\phi$-trivial 2-cells that intersect $P$, and
let $\partial D$ be those 1-cells that $D$ runs over.
The exact sequence associated with the pair $(D \cup \Gamma_n, \Gamma_n)$
contains the following terms:
$$H^1(D \cup \Gamma_n, \Gamma_n;{\Bbb F}_p) \rightarrow H^1(D \cup \Gamma_n;{\Bbb F}_p) \rightarrow
H^1(\Gamma_n;{\Bbb F}_p).$$
By excising $\Gamma_n - D$ from $(D \cup \Gamma_n, \Gamma_n)$, we see that
the first term is isomorphic to $H^1(D, \partial D; {\Bbb F}_p)$, which is
trivial. Hence,
$$H^1(D \cup \Gamma_n;{\Bbb F}_p) \rightarrow H^1(\Gamma_n;{\Bbb F}_p)$$ 
is an injection. So, the kernel of
$$H^1(\overline K_n; {\Bbb F}_p) \rightarrow 
H^1(\Gamma_n; {\Bbb F}_p)$$
equals the kernel of 
$$H^1(\overline K_n; {\Bbb F}_p) \rightarrow 
H^1(D \cup \Gamma_n; {\Bbb F}_p).$$
This clearly lies in the kernel of
$$H^1(\overline K_n; {\Bbb F}_p) \rightarrow 
H^1(A_n \cap B_n; {\Bbb F}_p)$$
because this latter map factors through $H^1(D \cup \Gamma_n; {\Bbb F}_p)$,
as $A_n \cap B_n$ is a subset of $D \cup \Gamma_n$.
This proves the claim.

The following claim verifies the key hypothesis of Theorem 1.3.

\noindent {\sl Claim 5.} The kernels of both the
following maps
$$\eqalign{
H^1(A_n; {\Bbb F}_p) & \rightarrow H^1(A_n \cap B_n; {\Bbb F}_p) \cr
H^1(B_n; {\Bbb F}_p) & \rightarrow H^1(A_n \cap B_n; {\Bbb F}_p) \cr}$$
have dimension at least 2, when $n$ is sufficiently large.

First note that $A_n \cap B_n$ has at most two components, since every edge
of $A_n \cap B_n$ is incident to $P$, which consists of two points.

The following is an extract of the Mayer-Vietoris sequence applied
to $A_n$ and $B_n$:
$$H^0(A_n \cap B_n; {\Bbb F}_p) \rightarrow H^1(\overline K_n; {\Bbb F}_p)
\rightarrow H^1(A_n; {\Bbb F}_p) \oplus H^1(B_n; {\Bbb F}_p).$$
Thus, the kernel of
$$H^1(\overline K_n; {\Bbb F}_p)\rightarrow H^1(A_n \cap B_n; {\Bbb F}_p)$$
maps onto a subspace of the direct sum of the kernels in Claim 5,
and the nullity of this mapping is at most the dimension of
$H^0(A_n \cap B_n; {\Bbb F}_p)$, which is at most 2. We therefore
know that the direct sum of the kernels in Claim 5 has dimension
at least $n(d(F) -1) - 1- c$, by Claim 4. So, to prove the
claim, it suffices to show that each of the kernels in Claim 5 has 
dimension at most $n(d(F) -1) - 3- c$. This will follow
from the claim below.

\noindent {\sl Claim 6.} The kernel of each of 
$$\eqalign{
H^1(A_n; {\Bbb F}_p) & \rightarrow H^1(A_n \cap B_n; {\Bbb F}_p) \cr
H^1(B_n; {\Bbb F}_p) & \rightarrow H^1(A_n \cap B_n; {\Bbb F}_p) \cr}$$
has dimension at most 
$\lceil n/2 \rceil (d(F)-1) + 2$.

We focus on the first of these kernels. The argument in the second
case is similar. Consider the exact sequence associated with the
pair $(A_n, A_n \cap B_n)$:
$$H^1(A_n, A_n \cap B_n; {\Bbb F}_p) \rightarrow H^1(A_n; {\Bbb F}_p)
\rightarrow H^1(A_n \cap B_n; {\Bbb F}_p).$$
Exactness implies that the kernel of the second map is equal to
the image of the first. It therefore suffices to bound the
dimension of $H^1(A_n, A_n \cap B_n; {\Bbb F}_p)$ from above.
Let $N(A_n \cap B_n)$ be a thin regular neighbourhood of $A_n \cap B_n$.
Then $H^1(A_n, A_n \cap B_n; {\Bbb F}_p)$ is isomorphic to $H^1(A_n, N(A_n \cap B_n); {\Bbb F}_p)$.
By excision, this is isomorphic to $H^1(A_n^-, \partial A_n^-)$,
where $A_n^- = A_n - {\rm int}(N(A_n \cap B_n))$. (See Figure 7
for an illustration of $A_n^- \cap C$ and $B_n^- \cap C$
within a 2-cell $C$.)

\vskip 18pt
\centerline{
\epsfbox{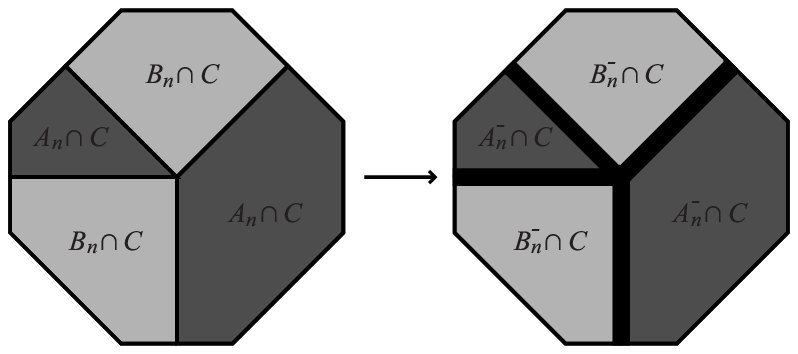}
}
\vskip 6pt
\centerline{Figure 7.}

The following
is an excerpt from the exact sequence for $(A_n^-, \partial A_n^-)$:
$$H^0(\partial A_n^-; {\Bbb F}_p) \rightarrow H^1(A_n^-, \partial A_n^-; {\Bbb F}_p)
\rightarrow H^1(A_n^-; {\Bbb F}_p).$$
This implies that the dimension of 
$H^1(A_n^-, \partial A_n^-; {\Bbb F}_p)$ is at most the
sum of the dimensions of $H^0(\partial A_n^-; {\Bbb F}_p)$
and $H^1(A_n^-; {\Bbb F}_p)$. Now, $\partial A_n^-$ has
at most two components, since each edge of $\partial A_n^-$ is
incident to one of two points in $X_n$ on the $A_n$ side of $P$.
Thus, to prove the claim, it suffices to show that
$$d_p(A_n^-) \leq \lceil n/2 \rceil (d(F) - 1).$$
The restriction of $A_n^-$ to any 2-cell $C$ that
intersects $P$ is a union of segments of $C$.
We may retract these segments onto their intersections with
the boundary of $C$, without
changing the homotopy type of $A_n^-$.
After this procedure, the 1-skeleton
of $A_n^-$ is $X_n \cap A_n^-$. Note that this is connected.
Thus,
$$d_p(A_n^-) \leq 1 - \chi(X_n \cap A_n^-)
= \lceil n/2 \rceil (d(F) - 1).$$
This proves the claim.

Claim 5 and Theorem 1.3 now give Theorem 1.2. $\square$

In the above proof, each of the 2-cells of
$\overline K_n$ was attached along a lift of a loop
$g_i^n$. For the $\phi$-non-trivial 2-cells,
it was important that the attaching word was
an $n^{\rm th}$ power. But, for the $\phi$-trivial 2-cells,
this was not important. It merely had to be of the form
$g_i^m$, for some $m \geq 2$. Then, working with cohomology
modulo $p$, for any prime $p$ dividing $m$, the argument works
unchanged. Also, we used a specific
homomorphism $\phi \colon F \rightarrow {\Bbb Z}$,
but it is not hard to see that any surjective
homomorphism would have worked. Thus, the same argument gives
the following variant of Theorem 1.2.

\noindent {\bf Theorem 2.1.} {\sl Let $F$ be a finitely
generated, free non-abelian group, and let $\phi \colon
F \rightarrow {\Bbb Z}$ be a surjective homomorphism.
Let $g_1, \dots, g_k$ be a collection of elements of
$F$ with trivial image under $\phi$, and let each of $g_{k+1},
\dots, g_r$ have non-trivial image. Then, there
is some integer $N$ with the following property.
For all $m \geq 2$ and $n \geq N$, the group
$F/ \langle \! \langle g_1^m, \dots, g_k^m,
g_{k+1}^n, \dots, g_r^n \rangle \! \rangle$
is large.}

\vskip 18pt
\centerline {\caps 3. Applications to Dehn surgery}
\vskip 6pt

It is a major conjecture that the fundamental group
of any finite-volume hyperbolic 3-manifold is large. This is
known to be true when the manifold has non-empty boundary by a theorem
of Cooper, Long and Reid (Theorem 1.3 of [1]). Many properties that
hold for 3-manifolds with toral boundary are also known
to be generically true for manifolds obtained by 
Dehn filling. Thus, `most' manifolds obtained by Dehn
filling a finite-volume hyperbolic 3-manifold should have
large fundamental group, but this remains conjectural
at present. However, we will show that if
it holds for one Dehn filling, then it holds for
infinitely many.

We now briefly recall the relevant terminology.
Let $M$ be a compact orientable 3-manifold with 
boundary a non-empty collection of tori. Then a
{\sl slope} on $\partial M$ is an isotopy class
of a simple closed curve which does not bound
a disc in $\partial M$. Let $(s_1, \dots, s_r)$
be a collection of slopes on $\partial M$,
with at most one $s_i$ on each component of $\partial M$.
Then $M(s_1, \dots, s_r)$ denotes the manifold
obtained by gluing on a collection of $r$ solid tori  to $M$
in the following way. We homeomorphically identify the boundary of
each solid torus $S^1 \times D^2$ with a boundary component of $M$,
so that $S^1 \times \partial D^2$ becomes a curve
of slope $s_i$. Then $M(s_1, \dots, s_r)$ is said to obtained 
from $M$ by {\sl Dehn filling}.

The operation of Dehn surgery is closely
related. Now one starts with a compact orientable
3-manifold $M$, which is possibly closed.
Pick a link $L$ in the interior of $M$, and
drill out a regular neighbourhood $N(L)$ of $L$,
which is a collection of solid tori. Then
Dehn filling $M - {\rm int}(N(L))$ along every component of
$\partial N(L)$ is called performing {\sl Dehn surgery}
on $L$.

For two slopes $s$ and $s'$ on a torus, their {\sl distance}
$\Delta(s,s')$ is the minimal intersection number between
representative simple closed curves.

\noindent {\bf Theorem 3.1.} {\sl Let $M$ be a compact orientable 
3-manifold with large fundamental group. Let $L$ be
a link in $M$, and let $\mu_1, \dots, \mu_r$ be a collection
of meridians for $L$, one for each component of the link. 
Then there is a positive integer $N$ with the following property.
Let $(s_1, \dots, s_r)$ be a collection of slopes,
one on each component of $\partial N(L)$, 
such that $N$ divides $\Delta(\mu_i, s_i)$ for each $i$. Then
the manifold obtained by performing Dehn surgery on $L$ via
the slopes $(s_1, \dots, s_r)$ has large fundamental group.
}

\noindent {\sl Proof.} Pick a basepoint for $M - {\rm int}(N(L))$, 
and, for each component of $\partial N(L)$, pick an arc
in $M - {\rm int}(N(L))$ joining it to the basepoint.
Use this to create, for each slope on $\partial N(L)$,
a fixed representative element of $\pi_1(M - {\rm int}(N(L)))$.
Let $\lambda_i$ be a curve on
$\partial M$ such that $\Delta(\mu_i, \lambda_i) = 1$.
Then, in the fundamental group of $M - {\rm int}(N(L))$, $s_i = \mu_i^{p_i}
\lambda_i^{q_i}$, where $|q_i| = \Delta(\mu_i, s_i)$.
Let $M'$ be the manifold obtained by
Dehn surgery along $L$ via the slopes $(s_1, \dots, s_r)$. Then 
$$\pi_1(M')/\langle \! \langle \mu_1, \dots, \mu_r
\rangle \! \rangle
\cong \pi_1(M)/\langle \! \langle \lambda_1^{\Delta(\mu_1, s_1)},
\dots, \lambda_r^{\Delta(\mu_r, s_r)} \rangle \! \rangle,$$
and this has as a quotient
$$\pi_1(M)/\langle \! \langle \lambda_1^N, \dots, \lambda_r^N
\rangle \! \rangle.$$
Now, for infinitely many $N$, the latter group is large,
by Theorem 1.1. Hence, for these values of $N$, so
is $\pi_1(M')$. $\square$

This has the following corollary.

\noindent {\bf Theorem 1.5.} {\sl Let $M$ be a compact orientable
3-manifold with boundary a collection of tori. Suppose that
for one collection of slopes $(s_1, \dots, s_r)$,
with one $s_i$ on each component of $\partial M$,
$M(s_1, \dots, s_r)$ has large fundamental group.
Then this is true for infinitely many distinct collections
$(s_1, \dots, s_r)$.}

\noindent {\sl Proof.} Suppose that for one collection of
slopes $(\mu_1, \dots, \mu_r)$, with one $\mu_i$ on each
component of $\partial M$, $M(\mu_1, \dots, \mu_r)$ has large
fundamental group. The cores of the filled-in solid tori
form a link $L$ in $M(\mu_1, \dots, \mu_r)$ with meridians
$\mu_1, \dots, \mu_r$. Apply Theorem 3.1 to obtain a positive integer $N$,
with the following property. For any collection of
slopes $(s_1, \dots, s_r)$,
with one $s_i$ on each component of $\partial M$,
$M(s_1, \dots, s_r)$ has large fundamental group,
provided that $\Delta(\mu_i,s_i)$ is a multiple of $N$
for each $i$. This clearly holds for infinitely many distinct collections
$(s_1, \dots, s_r)$. $\square$

\vskip 18pt
\centerline{\caps References}
\vskip 6pt

\item{1.} {\caps D. Cooper, D. Long, A. Reid,} 
{\sl Essential closed surfaces in bounded $3$-manifolds,} 
J. Amer. Math. Soc. 10 (1997) 553--563.

\item{2.} {\caps T. Delzant},
{\sl Sous-groupes distingu\'es et quotients des groupes hyperboliques.}
Duke Math. J. 83 (1996) 661--682.

\item{3.} {\caps N. Dunfield, W. Thurston,}
{\sl The virtual Haken conjecture: Experiments and examples,}
Geom. Topol. 7 (2003) 399--441.

\item{4.} {\caps M. Gromov}, {\sl Hyperbolic groups.}
Essays in group theory (1987) 75--263  Math. Sci. Res. Inst. Publ.

\item{5.} {\caps M. Lackenby,} {\sl Some 3-manifolds
and 3-orbifolds with large fundamental group}, Proc. Amer. Math. Soc. (to appear)

\item{6.} {\caps M. Lackenby,} {\sl 
 Large groups, Property $(\tau)$ and the homology growth of subgroups},
Preprint.

\vskip 12pt
\+ Mathematical Institute, University of Oxford, \cr
\+ 24-29 St Giles', Oxford OX1 3LB, United Kingdom. \cr

\end